\documentclass[12pt,draftcls,onecolumn]{IEEEtran}

\usepackage{graphics} 
\usepackage{epsfig} 
\usepackage{mathptmx} 
\usepackage{times} 
\usepackage{amsmath} 
\usepackage{amssymb}  

\newtheorem{thm}{Theorem}

\newcommand{\RR}{{\mathbb R}}

\newcommand{\HH}{{\mathbb H}}
\newcommand{\SSS}{{\mathbb  S}}
\newcommand{\CC}{{\mathbb C}}

\newcommand{\dotex}{{\frac{d}{dt}}}
\newcommand{\XX}{{\mathcal X}}
\newcommand{\YY}{{\mathcal Y}}

\title{\LARGE \bf Flatness-based control of a single qubit gate}

\author{
Paulo Sergio Pereira da Silva and Pierre Rouchon%
\thanks{The work of the first author was partially supported by CNPq. Both authors
were partially supported by CAPES/COFECUB.}
\thanks{P. S. Pereira da Silva is with Universidade de S\~ao
Paulo, Escola Polit\'{e}cnica - PTC -- Av. Prof. Luciano
Gualberto, Trav. 03, no. 158 05508-900 S\~AO PAULO - SP BRAZIL
{\tt \small paulo@lac.usp.br} }
\thanks{P. Rouchon is with Ecole des Mines de Paris, Centre Automatique et Syst\`{e}mes,
 60 Bd Saint-Michel, 75272 Paris cedex 06, FRANCE.
  {\tt \small pierre.rouchon@ensmp.fr}}%
}

\begin{document}
\maketitle

\begin{abstract}
This work considers the open-loop control problem of steering a
two level quantum system {from} an initial to a final condition.
The model of this system evolves on the state space $\mathcal{X}=
SU(2)$, having two inputs that correspond to the complex amplitude
of a resonant laser field. A symmetry preserving flat output is
constructed using a fully geometric construction and quaternion
computations. Simulation results
of this flatness-based open-loop control  are provided.\\

 \end{abstract}


\section{Introduction}
Take a single qubit, i.e. a two level quantum system. Denote by
$\omega_0$ its transition frequency. Assume that it is controlled
via a resonant laser field $v\in\RR$:
\begin{equation}\label{laser:eq}
    v= u\exp(-\imath \omega_0 t) + u^\ast \exp(\imath \omega_0 t)
\end{equation}
where $u=u_1+\imath u_2\in\CC$, $(u_1,u_2)\in\RR^2$, is its
complex amplitude. In general, the frequency  $\omega_0$ is  large
and the time variation of $u$ is slow: $|\dot u | \ll \omega_0
|u|$. In the interaction frame, after the  rotating wave
approximation and up to some scaling (see
e.g.,~\cite{haroche-raimond:book06}), the Hamiltonian reads $ u_1
\sigma_1 + u_2 \sigma_2$ where $\sigma_1$ and $\sigma_2$ are the
first two Pauli matrices (see appendix). The gate generation
problem then reads: take a transition  time $T>0$ such that
$\omega_0 T \ll 1$ and a goal matrix $\bar U\in SU(2)$; find a
smooth laser impulsion $[0,T]\ni t\mapsto u(t)\in\CC$ with
$u(0)=u(T)=0$ such that the solution $[0,T]\ni t \mapsto U(t)\in
SU(2)$  of the initial value problem
\begin{equation}\label{DynPauli:eq}
\imath \dotex U(t) = (u_1(t) \sigma_1 + u_2(t)\sigma_2)~ U(t),\quad
U(0)=I_2
\end{equation}
 reaches $\bar U$ at time $T$, i.e., $U(T)=\bar U$. This
motion planning problem admits a well known elementary
solution\footnote{The so-called ZYZ quantum logic gate.}. It
relies on the fact that $\bar U=$ $\exp(-\imath \gamma \sigma_1)$
$\exp(-\imath \beta \sigma_2)$ $\exp(-\imath \alpha \sigma_1)$,
for all $\bar U \in SU(2)$, for convenient
$(\alpha,\beta,\gamma)\in\RR^3$ (see,
e.g.,~\cite{altafini:qip02}). An obvious steering control $u(t)$
is decomposed into three elementary and successive pulses: for the
first (resp. third) pulse, $u_2=0$ and $u_1$ is such that its
integral over the pulse interval equals $\alpha$ (resp. $\gamma$);
for the second pulse, $u_1=0$ and the integral of $u_2$ is
$\beta$.

Here, we propose another solution  where $u_1$ and $u_2$ vary
simultaneously, i.e., the  steering control $u(t)$ is contained in
a single pulse. Our solution does not rely on optimal control
techniques (see for instance \cite{boscain-mason:jmp06} and the
references therein) and  is explicit. It does not rely on
numerical resolution scheme. It provides  controls that can be
chosen to be $C^\omega$ or $C^{\infty}$ function of $t$. As far as
we know, such explicit solution is new and could be of some
interest for reducing the transition time $T$ while still
respecting the rotating wave approximation.  Our approach is based
on the fact that the system dynamics is differentially
flat~\cite{fliess-et-al-ijc95}. The flat output constructed in
this paper has a clear geometrical interpretation.

In section~\ref{flatsym:sec}, theorem~\ref{flat:thm} shows, using
a quaternion description of~\eqref{DynPauli:eq}, that this system
is flat.
 We propose a coordinate free definition of
the flat-output that lives in the homogenous space
$SU(2)/\exp(\imath \RR \sigma_1)$. This geometric construction
preserves invariance with respect to right translations. In the
sense of~\cite{martin-et-al-cocv03}, the flat output is compatible
with right translations. The proposed construction can be seen as
the analogue of the geometric construction based on the Frenet
formula for the car system, where the steering angle is directly
related to the curvature of the path followed by the flat-output
curve \cite{rouchon-rudolph-ncn99}. In section~\ref{planing:sec},
we show how to use such geometric flatness  parameterization to
solve analytically the motion planning problem corresponding to an
arbitrary quantum gate. Simulations illustrate
theorem~\ref{steer:thm} and  the interest of such explicit
open-loop steering control. In section \ref{conclusions:sec}, some
conclusions are briefly stated. Some material has been deferred to
the appendix~\ref{pauli:ap}. In part A one finds the basics
properties of Pauli matrices and their associated quaternions as
well as the correspondence between $SU(2)$ and quaternions of
length one. In part B one finds a proof of the fact that the
motion planning algorithm has no singularities.

\section{A symmetry preserving flat output} \label{flatsym:sec}

The dynamics~\eqref{DynPauli:eq} reads in quaternion notation (see
appendix A)
\begin{equation}\label{dyn:eq}
 \dotex q = (u_1 e_1 + u_2 e _2) q
\end{equation}
where $q\in \HH_1$ is a quaternion of  length one and where
$(u_1,u_2) \in\RR^2$ is the control relative to the modulation of a
coherent laser field ($u_1+\imath u_2 $ is the complex field
amplitude). This system is a driftless system on the Lie Group
$\HH_1$. It is
 controllable (see, e.g.,
\cite{alessandro1}). Moreover, this control system is invariant
with respect to right translations in the sense
of~\cite{martin-et-al-cocv03,bonnabel-et-al:acc06}:
\begin{itemize}
   \item the group $G=\HH_1$ acts on the state space $\XX=\HH_1$ via
right multiplication $\phi_g: q\mapsto q g$ where $q\in \HH_1$.
   \item the dynamics is $G$-invariant: if $t\mapsto
   (q(t),u_1(t),u_2(t))$ is a solution of~\eqref{dyn:eq} then
   $t\mapsto (q(t) g,$ $u_1(t),$ $u_2(t) )$ is also a solution
   of~\eqref{dyn:eq} for any $g\in G$.
\end{itemize}
The controllability structure of this system is in fact  of a very
special kind. Around any point $\bar q\in \HH_1$, \eqref{dyn:eq} can
be seen in local coordinates as a driftless controllable system with
$3$ states\footnote{Take e.g., the exponential map: $
 (x_1,x_2,x_3)\in\RR^3 \mapsto
  \exp(x_1 e_1+x_2 e_2+x_3 e_3) \bar q
$ that maps  a neighborhood of $0\in\RR^3$ to a neighborhood of
$\bar q$ in $\HH_1$.} and 2 controls. Thus, as known
since~\cite{charlet-et-al-siam91} (see
also~\cite{martin-et-al-caltech03}),  such system is
differentially flat and the flat output function can be chosen to
depend only on the state. More precisely,  the  flat output for
the controllable system $\dotex x= u_1 f_1(x)+ u_2 f_2(x)$ with
$\dim(x)=3$ is obtained by the rectifying coordinates of any
vector field $f(x)= \alpha_1(x) f_1(x) + \alpha_2(x) f_2(x) $
which is a linear combination of the two control vector fields
$f_1$ and $f_2$ ($\alpha_1$, $\alpha_2$ are any scalar functions
of $x$).

We  propose here  a coordinate free and symmetry preserving
construction of the flat output via the previous procedure. Thus
we are looking for a flat output map $h:\HH_1 \mapsto \YY$, where
$\YY$ is the output space, a compact manifold of  dimension $2$,
and $G$-compatible in the sense of~\cite{martin-et-al-cocv03}.
This means that the output map $h$ must satisfy the following
constraint: there exists an action of $G=\HH_1$ on the flat output
space $\YY$ described by the transformation group $\varrho_g:
y\mapsto \varrho_g (y)$ such that $\varrho_g(h(q))=h(qg)$ for any
$q\in \HH_1$.  The following  construction will be based on the
control vector field associated to $u_1$, and hence to $e_1$.

Denote by $K=\{\exp(\phi e_1)\}_{\phi\in[0,2\pi]}$ the one
dimensional subgroup of $\HH_1$ generated  by $e_1$. We can
consider the action of $K$ on $\HH_1$ via left multiplication: to
any $k\in K$, we have the diffeomorphism $\HH_1\ni q \mapsto kq\in
\HH_1$. Two elements of $\HH_1$, $q$ and $p$,  belong to the same
orbit if and only if there exists $k\in K$ such that $kq=p$.
Denote by $\YY$ the set of the orbits. This set is a compact
manifold of dimension 2 and the output function $h$ is the map
that associates to any $q$, the orbit $h(q)$ to which $q$ belongs.
This map is a smooth submersion, and $\YY$ is called an homogenous
space (see, e.g.,~\cite{olver-book95}).  If $q$ and $p$ belong to
the same orbit, $qg$ and $pg$ also belong to the same orbit for
any $g\in \HH_1$. Therefore, this output map is $G$-compatible in
the sense of~\cite{martin-et-al-cocv03}.

Assume that $y(t)$ is a curve on $\mathcal Y = \HH_1 / K$, at
least of class $C^2$. Since the map $h : \HH_1 \rightarrow
\mathcal Y$ is a submersion, in adequate local coordinates one has
$h(x_1,x_2,x_3) = (x_1, x_2)$. Assume, without loss of generality,
that the open neighborhood of definition of h is
\emph{rectangular} and contains $(0,0,0)$. Define locally the
smooth map $g : U \subset \mathcal Y \rightarrow V \subset \HH_1$,
where $U$, $V$ are open sets and $g(x_1, x_2) = (x_1, x_2, 0)$.
Note that $g$ is smooth, and $Y(t)= g(y(t))$ is such that $h(Y(t))
= y(t)$. Then, locally, there exist smooth maps $g^{(1)}$ and
$g^{(2)}$ such that $\dot{Y}(t) = g^{(1)}(y(t), \dot y(t))$ and
$\ddot{Y}(t) = g^{(2)}(y(t), \dot y(t), \ddot y(t))$.

Let us show now that the map $h$ defines a flat output. This means
that the inverse of  system~\eqref{dyn:eq} with  output $y=h(q)$
has no dynamics\footnote{This is equivalent to say that the state
$q$ and the input $u =(u_1, u_2)$ can be written respectively as
$q={\cal A}(y,\dot y, \ddot y, \ldots, y^{(\alpha)})$ and $u =
{\cal B}(y,\dot y, \ddot y, \ldots, y^{(\beta)})$ for convenient
smooth maps $\cal A$ and $\cal B$.}. Thus we have to consider the
following implicit system
$$
 \dotex  q = (u_1 e_1 + u_2 e_2) q, \quad y = h(q)
$$
where $t\mapsto y(t) $ is a known function of time and where the
quaternion $q(t)\in \HH_1$ and the control $(u_1(t),u_2(t))$ are
the unknown quantities.

The problem is how to manipulate $h$, since  only a geometric
construction for $h$ is available. Knowing the function $t\mapsto
y(t)$ means that we have at our disposal a smooth  function
$t\mapsto Y(t)\in \HH_1$  such that $y(t)=h(Y(t))$. Hence, to have
$y(t)=h(q(t))$ means that $q$ and $Y$ belongs to the same orbit
for each time $t$. Therefore, there exists $k(t)=\exp(
\phi(t)e_1)$ in $K$ such that $q = k  Y$.  Since $k(t) = q(t)
Y^*(t)$, then $k(t)$ is smooth. Thus, we have
$$
\dotex q = \left(\dotex k\right) Y + k \dotex Y .
$$
But $\dotex k =\omega e_1 k$ where $\omega=\dotex \phi$.
Using~\eqref{dyn:eq}, we get the following equation $ k\dotex Y=
((u_1-\omega) e_1 + u_2 e_2) k Y $, that is
$$
  k\left(\dotex Y\right)Y^\ast k^\ast= (u_1-\omega)e_1 + u_2 e_2
  .
$$
 This quaternion equation
gives in fact  $k$ as a function of $\left(\dotex Y\right) Y^\ast$.
Left and right multiplication by $e_3$ yields
$$
e_3 k \left(\dotex Y\right) Y^\ast k^\ast e_3 = (u_1-\omega) e_1 +
u_2e_2
$$
since $e_3 e_i e_3 =e_i$ for $i=1,2$. Hence, we have the following
relation (without the controls and $\omega$):
\begin{equation}\label{ydot0:eq}
e_3 k \left(\dotex Y\right) Y^\ast k^\ast e_3 =   k \left(\dotex
Y\right) Y^\ast k^\ast.
\end{equation}
Assume that
\begin{equation}\label{ydot:eq}
 \left(\dotex Y\right) Y^\ast =  \omega_1 e_1 + \omega_2 e_2
 + \omega_3 e_3
\end{equation}
 where the $\omega_i$'s are known smooth
real functions of time. Thus, we get
$$
k \left(\dotex Y\right) Y^\ast k^\ast =\omega_1 e_1 + k^2
(\omega_2e_2+\omega_3 e_3)
$$
since $e_1 k^\ast = k^\ast e_1$, $kk^\ast=1$ and $e_ik^\ast= k
e_i$ for $i=2,3$.  Therefore, \eqref{ydot0:eq} reads:
$$
k^4 (\omega_2 e_2 + \omega_3 e_3) = (\omega_2 e_2 - e_3 e_3)
$$
since $e_3 k^2=(k^\ast)^2 e_3$ and $k^{-1}=k^\ast$.

Right multiplication by $e_2$  yields the following algebraic
equation defining~$k$
$$
k^{4} (\omega_2  +  \omega_3 e_1) = (\omega_2  - \omega_3 e_1).
$$
Since  $k=\cos\phi + \sin\phi e_1$, we have the following equation
for the angle $\phi$
$$
(\cos 4\phi +  \sin 4\phi ~e_1 ) (\omega_2  + \omega_3 e_1) =
(\omega_2  -\omega_3 e_1)
$$
which is equivalent to $ \exp( 4\phi\imath )= \frac{z^2}{|z|^2}$
where $z=\omega_2 - \omega_3 \imath $ is a known complex number.
Thus, there are four distinct possibilities  for  $k$:
\begin{equation}\label{Uyydot:eq}
k=\pm \exp\left(\frac{\theta}{2} e_1\right), \quad k= \pm e_1
\exp\left(\frac{\theta}{2} e_1\right)
\end{equation}
where $\theta$ is the argument of  $\omega_2-\omega_3\imath $. The
controls $u_1$ and $u_2$ associated to one of these four
trajectories are obtained by
$$
  e_3 k \dotex Y Y^\ast k^\ast e_3 =
  (u_1-\omega) e_1 + u_2e_2
$$
where $2\omega= \dotex \theta$ is given via simple algebraic
formulae based on $\omega_2$, $\omega_3$, $\dotex \omega_2$ and
$\dotex \omega_3$:
$$
 \omega = \frac{\omega_3\dotex \omega_2-\omega_2 \dotex \omega_3 }
   {2(\omega_2^2+\omega_3^2)}
   .
$$
For the two  branches  $k=\pm \exp\left(\frac{\theta}{2} e_1\right)$
we get
$$
 \left\{
 \begin{aligned}
 u_1&= \omega_1 +
    \frac{\omega_3\dotex \omega_2-\omega_2 \dotex \omega_3 }
    {2(\omega_2^2+\omega_3^2)}
  \\
  u_2&= \sqrt{\omega_2^2+\omega_3^2} \end{aligned}
 \right.
$$
and for the two other ones $k=\pm e_1 \exp\left(\frac{\theta}{2}
e_1\right)$ we get
$$
 \left\{
 \begin{aligned}
 u_1&= \omega_1 +
    \frac{\omega_3\dotex \omega_2-\omega_2 \dotex \omega_3 }
    {2(\omega_2^2+\omega_3^2)}
  \\
  u_2&= -\sqrt{\omega_2^2+\omega_3^2} \end{aligned}
 \right.
$$
where  just the sign of $u_2$ is changed. All the previous
computations are valid when $\omega_2-\omega_3\imath \neq 0$,
i.e., when $\dotex y \neq 0$: $(\omega_2^2+\omega_3^2)$  does not
depends on $Y(t)$ such that $h(Y(t))=y(t)$; it depends only on
$y(t)$ and vanishes if, and only if, $\dotex y(t)=0$. To
summarize, we have proved the following result:
\begin{thm} \label{flat:thm}
Take $T>0$ and  an arbitrary $C^2$ curve $[0,T]\ni t\mapsto y(t)$ on
$\YY$ such that $\dotex y(t)\neq 0$ for any $t\in[0,T]$. For any
smooth curve  $t\mapsto Y(t)\in\HH_1$ such that $h(Y(t))=y(t)$,  set
 $z=\omega_2(t)-\omega_3(t)\imath \neq 0$ for all $t\in[0,T]$
where  $\left(\dotex Y\right)Y^\ast=\omega_1 e_1+\omega_2
e_2+\omega_3 e_3$. Then there exists a smooth function $[0,T]\ni
t\mapsto \theta(t)\in\RR$   such that $\exp( \theta\imath )=
\frac{z}{|z|}$ and any smooth solution $t\mapsto
(q(t),u_1(t),u_2(t))$ of~\eqref{dyn:eq} satisfying $h(q(t))=y(t)$
for all $t\in[0,T]$ is one of the four following trajectories
indexed by $n\in\{0,1,2,3\}$:
\begin{equation}\label{flat:eq}
 \left\{
 \begin{aligned}
   q(t)& = (e_1)^n \exp\left(\frac{\theta(t)}{2}~e_1\right) Y(t)
  \\
  u_1&= \omega_1 +
    \frac{\omega_3\dotex \omega_2-\omega_2 \dotex \omega_3 }
    {2(\omega_2^2+\omega_3^2)}
  \\
  u_2&= (-1)^n \sqrt{\omega_2^2+\omega_3^2}
 \end{aligned}
 \right.
\end{equation}
\end{thm}

Recall that some $Y(t)$ such that $y(t)=h(Y(t))$ is locally given
by $Y(t) = g(y(t))$,  and furthermore $\dot Y(t)=g^{(1)}(y, \dot
y)$ and $\ddot Y(t)=g^{(2)}(y, \dot y, \ddot y)$, where $g$,
$g^{(1)}$ and $g^{(2)}$ are smooth maps. In particular, the last
theorem proves that $y=h(q)$ is a flat output.

The flat output $y=h(q)$ is obtained with $e_1$ playing a specific
role. In fact one can see that any map $h_\eta:\HH_1 \mapsto
\HH_1/K_\eta$ ($\eta\in\SSS^1$) corresponding to the subgroup
$K_\eta=\exp(\RR (\cos\eta e_1 +\sin\eta e_2))$  also defines a
flat output. It just corresponds to a rotation by the angle $\eta$
of $(q_1,q_2)$ and $(u_1,u_2)$. If we set
$$
 e_1=\cos\eta \tilde e_1 + \sin \eta \tilde e_2, \quad
 e_2= -\sin\eta \tilde e_1 + \cos \eta \tilde e_2
$$
the imaginary quaternions $(e_1,e_2,e_3)$ and $(\tilde e_1,\tilde
e_2,e_3)$ satisfy exactly the same commutation relations. Thus, if
$t\mapsto q(t)$ is a solution of~\eqref{dyn:eq} with the control
$(u_1(t),u_2(t))$ then
\[
     t\mapsto q_0(t) + (\cos\eta q_1(t) -
\sin\eta q_2(t))e_1 + (\sin\eta q_1(t) + \cos\eta q_2(t)) e_2 +
q_3(t) e_3
\]
is also a solution of~\eqref{dyn:eq} with the control
$$
 \tilde u_1 = \cos\eta u_1(t) - \sin\eta u_2(t)
  ,\quad
 \tilde u_2 = \sin\eta u_1(t) + \cos\eta u_2(t)
 .
$$
This symmetry and the fact that, as stated in
theorem~\ref{flat:thm}, $h=h_0$ is a flat output, implies directly
that $h_\eta$ is also a flat-output. The family
$(h_\eta)_{\eta\in\SSS^1}$  is made of flat outputs  all compatible
versus right translations.

\section{Motion planning} \label{planing:sec}
In this section,  we will use~\eqref{flat:eq} with $n=0$ to
propose an explicit solution for the motion planning problem
stated in the introduction: for any $T>0$ and any final state
$\bar q\in\HH_1$, find a smooth control $[0,T]\ni t \mapsto
u(t)=(u_1(t),u_2(t))\in\RR^2$ with $u(0)=u(T)=0$, such that the
solution $[0,T]\ni t \mapsto q(t)\in\HH_1$ of~\eqref{dyn:eq}
starting {from}  $q(0)$ reaches $\bar q$ at time $T$: i.e.,
$q(T)=\bar q$.

 As the system is driftless, every
time re-parameterization of a solution is also a solution. In fact,
consider the equation
\[
\frac{d}{ds} \tilde q(s) = (\tilde {u}_1(s) e_1 + \tilde
 {u}_2(s) e_2) \tilde q(s)
\]
Let $\varsigma : [0, T] \rightarrow [0, 1]$ be an increasing
diffeomorphism. Then $\tilde q(s)$ is a solution of the previous
equation defined on $[0, 1]$, with input $(\tilde {u}_1(s) ,
\tilde {u}_2(s))$ if and only if $q(t) = \tilde q (\varsigma(t))$
is a solution of (\ref{dyn:eq}) defined on $[0, T]$ with input
$(u_1(t), u_2(t))= \frac{d \varsigma}{dt} ( \tilde {u}_1
(\varsigma(t)) , \tilde {u}_2 (\varsigma(t))$. One concludes that,
without loss of generality, one may always state the motion
planning problem with the (virtual) time $s$ belonging to the
interval $[0,1]$ and after that, one may ``control the clock'' by
choosing a convenient bijection $s=\varsigma(t)$. Thus, it is
enough to solve the motion planning problem in the $s$ scale where
we can disregard  the fact that the control has to vanish  at the
beginning and at the end: it is enough to take for example
$\varsigma(t)= 3 \left(\frac{t}{T}\right)^2-2
\left(\frac{t}{T}\right)^3$ to get $u$ equal to zero at $t=0$ and
at $t=T$, since $\dotex \varsigma(0)=\dotex \varsigma(T)=0$.

In the sequel we propose a solution in the $s$-scale. For
clarity's sake,  we will remove the $\tilde~$ when $u$ and $q$ are
considered as function of $s$. The derivation in $s$ will be
denoted by $^\prime$: $du/ds=u^\prime$, $dq/ds=q^\prime$, \ldots

Thus, we have to find a smooth control $[0,1]\ni s \mapsto u(s)$
such that the solution of
$$
 q^\prime(s) = (u_1(s) e_1 + u_2(s) e_2) q(s), \quad q(0)=1
$$
satisfies $q(1)=\bar q$, where $\bar q$ is any goal state in
$\HH_1$.

We can always assume that
$$
\bar q = \bar q_0 + \sqrt{\bar q_1^2+\bar q_2^2} (\sin\bar\eta
e_1+\cos\bar\eta e_2) + \bar q_3
$$
for some angle $\bar\eta\in [0,2\pi]$. Thus, as explained at the
end of last section, up to a rotation of angle $\bar\eta$ of the
control, we can assume  that $\bar q_1=0$. More precisely, if
$\bar q_1\neq 0$, set $\bar\eta$ to be the argument of the complex
$\bar q_2 + \bar q_1\imath $. If $s\mapsto (u_1(s),u_2(s))$ steers
$q$ {from} $q(0)=1$ to $q(1)=\bar q_0 +  \sqrt{\bar q_1^2+\bar
q_2^2} e_2 + \bar q_3 e_3$, then the control
$$
 s\mapsto
  (\cos\bar\eta u_1(s) + \sin\bar\eta u_2(s),
   -\sin\bar\eta u_1(s) + \cos\bar\eta u_2(s) )
$$
steers $q$ {from} $q(0)=1$ to $q(1)=\bar q$.

Thus up-to a rotation of angle $\bar\eta$ of the control, we can
assume that $\bar q_1=0$ and $\bar q_2 \geq 0$. Thus we can define
two angles $\bar\alpha\in]0,\pi]$ and $\bar\beta\in
[-\frac{\pi}{2},\frac{\pi}{2}]$ such that
$$
\bar q = \cos\bar\alpha + \sin\bar\alpha (\cos\bar\beta
e_2+\sin\bar\beta e_3).
$$
If the control $s\mapsto u(s)$ steers the system {from} $q(0)=1$
to $q(1)=\cos\bar\alpha + \sin\bar\alpha (\cos\bar\beta
e_2+\sin\bar\beta e_3)$, the same control steers the system {from}
$$
q(0)=\cos\bar\lambda + \sin\bar\lambda (\cos\bar\beta
e_2+\sin\bar\beta e_3)
$$
to
$$
q(1)=\cos(\bar\lambda+\bar\alpha) + \sin(\bar\lambda+\bar\alpha)
(\cos\bar\beta e_2+\sin\bar\beta e_3) .
$$
This is a direct consequence of right translation invariance
of~\eqref{dyn:eq} and right multiplication by $\cos\bar\lambda +
\sin\bar\lambda (\cos\bar\beta e_2+\sin\bar\beta e_3)$.

Take now the formulae~\eqref{flat:eq} in the $s$-scale with
\begin{equation}
 \label{flat:eq2}
Y(s)= \cos (\alpha(s)) + \sin (\alpha(s)) (\cos(\beta(s))e_2 +
\sin(\beta(s)) e_3
\end{equation}
where $\alpha(s)$ and $\beta(s)$ are smooth functions such that
\begin{equation}\label{ab:eq}
 \alpha(0)=\bar\lambda,\quad \alpha(1)=\bar\lambda+\bar\alpha
 ,\quad \beta(0)=\beta(1)=\bar\beta
 .
\end{equation}
Set, as in theorem~\ref{flat:thm}
$$
Y^\prime Y^\ast = \omega_1(s) e_1 + \omega_2(s) e_2 +\omega_3(s) e_3
.
$$
Simple computations shows that
$$
z = \omega_2 -\omega_3\imath = \exp(-\imath\beta)(\alpha^\prime -
\imath \beta^\prime \cos\alpha\sin\alpha) .
$$

 Now we shall construct \eqref{flat:eq2}
 such that $q(s)= \exp(\phi(s) e_1) Y(s), s \in [0, 1]$
 is a trajectory of
 the system. We will assume that $q(0)= Y(0)$ and $q(1)=Y(1)$.
 So we must have $\phi(0) = \phi(1) =0$. Furthermore, if we can
ensure that $s\mapsto z(s)$ never vanishes, and
$\theta(0)=\theta(1)=0$, then the trajectory of~\eqref{flat:eq}
with $n=0$ will provide a steering control $u$.

Let us now show in detail how to design  the functions $\alpha(s)$
and $\beta(s)$ satisfying these constraints. First of all we have
the initial and final constraints~\eqref{ab:eq}. By taking
$$
\bar\lambda = \left\{
           \begin{array}{ll}
              - \frac{\bar\alpha}{2}, & \hbox{for
                         $\bar\alpha\in[\frac{\pi}{4},\frac{3\pi}{4} ]$;} \\
             \frac{\pi}{4} - \frac{\bar\alpha}{2} , & \hbox{otherwise;}
           \end{array}
         \right.
$$
we always have $\cos\alpha\sin\alpha$ far {from} $0$ when $s=0$
and $s=1$. Thus we can impose the following  initial and final
constraints  for $\beta^\prime$:
$$
 \beta^\prime(0)=-\frac{\bar\alpha\sin\bar\beta}
                       {\sin{\bar\lambda} \cos{\bar\lambda}}
 \; \quad
 \beta^\prime(1)=-\frac{\bar\alpha\sin\bar\beta}
                       {\sin({\bar\lambda}+\bar\alpha)\cos({\bar\lambda}+\bar\alpha)}
$$
and for $\alpha^\prime$
$$
\alpha^\prime(0)=\alpha^\prime(1)=\cos\bar\beta \bar\alpha .
$$
Then $\alpha(s)$ and  $\beta(s)$  are the  polynomials of degree
$\leq 3$ satisfying these initial and final  constraints. Since
$\bar\alpha >0$ and $|\bar\beta|\leq \frac{\pi}{2}$,
$s\mapsto\alpha(s)$ can be a strictly increasing function on
$[0,1]$ and
 $\alpha^\prime >0$ for $s\in]0,1[$ (see appendix B). Thus  the complex number
$$
z=\exp(-\imath\beta)(\alpha^\prime - \imath \beta^\prime
\cos\alpha\sin\alpha)
$$
never  vanishes for $s\in]0,1[$.   For $s=0$ and $s=1$, we have
$$
\alpha^\prime - \imath \beta^\prime \cos\alpha\sin\alpha =
\exp(\imath \bar\beta) \bar\alpha .
$$
Thus $z(0)=z(1)=\bar\alpha > 0$. To summarize the closed path
$[0,1]\ni s\mapsto z(s)\in\CC $ never passes through $0$ nor turns
around $0$. We satisfy the assumption of theorem~\ref{flat:thm} in
the $s$-scale. Moreover we can set
$z(s)=r(s)\exp(\imath\theta(s))$ with $r(s) >0$ and $\theta(s)$
smooth functions on $[0,1]$ with $\theta(0)=\theta(1)=0$. We avoid
with such design of $\alpha(s)$ and $\beta(s)$ the monodromy
problem associated to  the  resolution of $(\exp(\imath \phi))^4 =
z^2/|z|^2$.  Finally we have proved the following result.

\begin{thm} \label{steer:thm}
Take $\bar q =\bar q_0 + \bar q_1 e_1 + \bar q_2 e_2 +\bar q_3
e_3\in\HH_1$ with $\bar q \neq  1$. Chose $\bar \eta\in[0,2\pi[$
such that $q_1 e_1 + q_2 e_2 = \sqrt{\bar q_1^2+\bar q_2^2}
(\sin\bar\eta e_1 + \cos\bar\eta e_2)$. Define
$\bar\alpha\in]0,\pi]$ and
$\bar\beta\in[-\frac{\pi}{2},\frac{\pi}{2}]$ such that
$$
\bar q_0 + \sqrt{\bar q_1^2+\bar q_2^2} ~e_2 + \bar q_3 e_3=
\cos\bar\alpha + \sin\bar\alpha (\cos\bar\beta e_2 + \sin\bar\beta
e_3) .
$$
Set $ \bar\lambda = - \frac{\bar\alpha}{2}$ if
                         $\bar\alpha\in[\frac{\pi}{4},\frac{3\pi}{4}
]$ and $\bar\lambda = \frac{\pi}{4} - \frac{\bar\alpha}{2}$
otherwise. Define $\alpha(s)$ and $\beta(s))$ as being the unique
polynomial functions of degree $\leq 3$ such that ($^\prime$ stands
for $d/ds$)
\begin{align*}
&\alpha(0)=\bar\lambda, \quad \alpha(1)=\bar\lambda+ \bar\alpha,
\quad \alpha^\prime(0)=\alpha^\prime(1)=\bar\alpha \cos\bar\beta
\\
&\beta(0)=\beta(1)=\bar\beta
\\
& \beta^\prime(0)=-\frac{\bar\alpha\sin\bar\beta}
                       {\sin\bar\lambda\cos\bar\lambda}
 ,\quad
 \beta^\prime(1)=-\frac{\bar\alpha\sin\bar\beta}
                       {\sin(\bar\lambda+\bar\alpha)\cos(\bar\lambda+\bar\alpha)}
\end{align*}
Define  $\omega_1(s)$, $\omega_2(s)$ and $\omega_3(s)$  by
\begin{align*}
\omega_1 &= (1 - 2 \cos^2(\alpha)) \beta^\prime
\\
 \omega_2 - \imath \omega_3 &= \exp(-\imath \beta)
(\alpha^\prime - \imath \beta^\prime \sin\alpha\cos\alpha) .
\end{align*}
 Then $\omega_2$ and $\omega_3$ never vanish
simultaneously  and the control
\[
  \left(
    \begin{array}{c}
     u_1(t)
      \\
     u_2(t)
    \end{array}
  \right)
   =
\dotex\varsigma(t)
\left(
  \begin{array}{cc}
    \cos\eta & \sin\eta\\
    -\sin\eta & \cos\eta
  \end{array}
\right) \left(
  \begin{array}{c}
    \omega_1 +
    \frac{\omega_3\omega_2^\prime-\omega_2\omega_3\prime }
    {2(\omega_2^2+\omega_3^2)}
     \\
  \sqrt{\omega_2^2+\omega_3^2}
  \end{array}
\right)_{s=\varsigma(t)}
\]
steers  system~\eqref{dyn:eq} {from} $q(0)=1$ to $q(T)=\bar q$
with $t \mapsto \varsigma(t)\in[0,1]$ being a $C^k$ increasing
bijection between $[0,T]$ and $[0,1]$ $k\geq 1$. When in addition
$\frac{d^n\varsigma}{dt^n}|_s=0$ for $s=0$ and $s=1$, and
$n=1,\ldots, k$, the control $t\mapsto u(t)$ is $C^{k-1}$ with
$\frac{d^{n-1}u}{dt^{n-1}}=0$ for $s=0$ and $s=1$.
\end{thm}

\begin{figure}
  \centerline{\includegraphics[width=.8\textwidth]{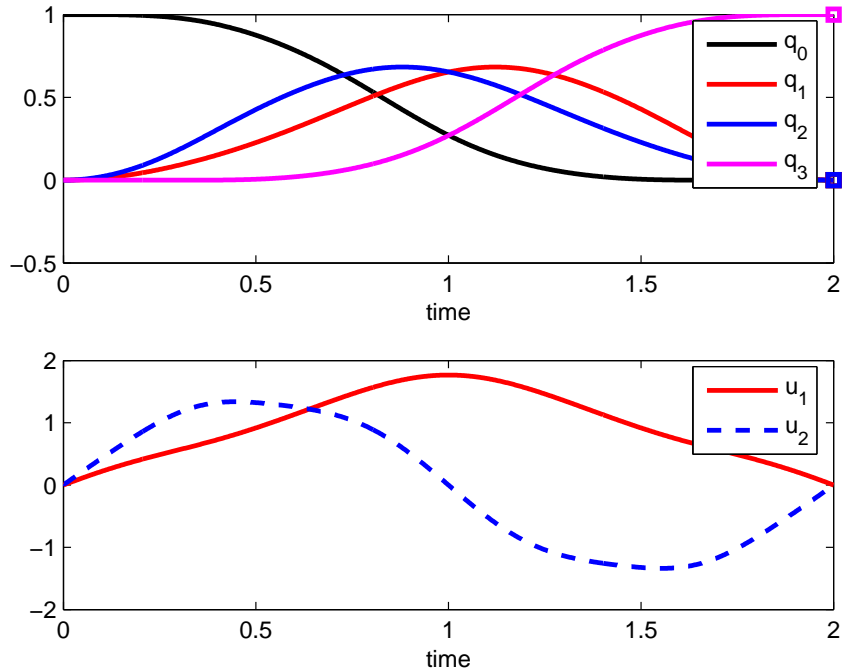}}
  \caption{The steering control and trajectory derived {from}
theorem~\ref{steer:thm} with $T=2$, $\bar q=e_3$ and $\varsigma(t)=3
(t/T)2 - 2 (t/T)3$. The control magnitude is very close to  an $ZYZ$
control design with two separated $\frac{\pi}{2}$ pulses.  The
simulation code  (matlab m-file  and scilab sci-file) can be
downloaded {from} {\tt
http://cas.ensmp.fr/~rouchon/publications/PR2007/CodeMatlabScilabQubit.zip}.
}\label{steer:fig}
\end{figure}

Figure~\ref{steer:fig}  illustrates the steering control described
by theorem~\ref{steer:thm} with $T=2$, $\bar q_0= e_3$, and
$\varsigma(t)= 3(t/T)^2-2(t/T)^3$.  We see that the control is a
smooth function  with maxima around $\pi/2$, a value close to the
$ZYZ$ design based on two successive pulses:
$(u_1,u_2)=(0,\frac{\pi}{2})$ for $t\in[0,1]$ and
$(u_1,u_2)=(\frac{\pi}{2}, 0)$ for $t\in[1,2]$. Thus our flatness
based design  yields, with the  same transition time and control
magnitude,  smooth control actions.

\section{Concluding remarks}
\label{conclusions:sec} The results of this paper holds if the
laser matches  exactly the resonant frequency. If we have a
frequency offset of $\Delta_r$ {from} resonance, then this offset
leads to the following drift (see, e.g.,
\cite{haroche-raimond:book06}):
$$
\dotex q = (u_1 e_1 + u_2 e_2 + \Delta_r e_3) q .
$$
It is still interesting to notice that  $h(q)$ is also a flat
output. In this case, the key relation~\eqref{ydot0:eq} becomes
$$
e_3 k \left(\dotex Y\right) Y^\ast k^\ast e_3 =   k \left(\dotex
Y\right) Y^\ast k^\ast + 2\Delta_r e_3.
$$
and $k=\exp(\phi e_1)$  is a root of the following polynomial
$$
k^4 (\omega_2 e_2 + \omega_3 e_3) + 2 k^2 \Delta_r e_3 - (\omega_2
e_2 - \omega_3 e_3) = 0 .
$$
Then one could try to apply similar techniques for solving the
motion planning problem for this system, although the time-scale
$s = \zeta(t)$ cannot be considered in this case.

\bibliographystyle{IEEEtran}

\appendix

\noindent\emph{A -- Pauli Matrices and Quaternions}

\label{pauli:ap} The Hermitian matrices
  $$
  \sigma_1=\left(
             \begin{array}{cc}
               0 & 1 \\
               1 & 0 \\
             \end{array}
           \right)
           ,\quad
  \sigma_2=\left(
             \begin{array}{cc}
               0 & -\imath \\
               \imath & 0 \\
             \end{array}
           \right)
           ,\quad
  \sigma_3=\left(
             \begin{array}{cc}
               1 & 0 \\
               0 & -1 \\
             \end{array}
           \right)
  $$
are the three Pauli matrices.  They satisfy  $\sigma_k^2=1$,
$\sigma_k \sigma_j = - \sigma_j \sigma_k$ for $k\neq j$, and
  $$
  \sigma_1\sigma_2 = \imath \sigma_3, \quad
  \sigma_2\sigma_3 = \imath \sigma_1, \quad
  \sigma_3\sigma_1 = \imath \sigma_2 .
  $$
  Any matrix $U$ in $SU(2)$ reads
  $$
  U = q_0 - q_1 \imath \sigma_1 - q_2 \imath\sigma_2 - q_3 \imath
  \sigma_3
  $$
  with $(q_0,q_1,q_2,q_3)\in\RR^4$ such that  $q_0^2+q_1^2+q_2^2+q_3^2=1$.
By setting
$$
e_1= -\imath \sigma_1, \quad e_2=-\imath \sigma_2, \quad e_3=-\imath
\sigma_3
$$
on can identify $SU(2)$ with the set of quaternions
$$
q=q_0+q_1 e_1 + q_2 e_2 + q_3 e_3
$$
of length one. This set is denoted by  $\HH_1$ and corresponds to
quaternions $q\in\HH$ such that $qq^\ast=1$ where  $q^\ast =q_0-q_1
e_1 - q_2 e_2 - q_3 e_3$ is the conjugate quaternion of $q$. Thus
the dynamics~\eqref{DynPauli:eq} becomes~\eqref{dyn:eq} with $q$
corresponding to $U$. Notice that $\HH_1$ is a compact Lie group of
dimension $3$.

We recall here some useful  relations for $k=1, 2, 3$, $j\neq k$
and $\phi\in\RR$:
\begin{align*}
& e_k^2=-1, \quad e_ke_j=-e_je_k,\quad  \exp(\phi e_k)=\cos\phi +
e_k\sin\phi
\\
&\exp(\phi e_k) e_j = e_j \exp(-\phi e_k)
\\
&e_1e_2=e_3, \quad e_2e_3=e_1, \quad e_3e_1=e_2
\end{align*}

\noindent\emph{B -- Proof that $z= \omega_2 - \imath \omega_3$
never vanishes for $s \in ]0,1[$}

Since $\omega_2 - \imath \omega_3  = \exp(\imath \beta)
(\alpha^\prime - \imath \beta^\prime \sin \alpha \cos \alpha)$,
 it suffices to show that $\alpha^\prime>0$ for $s \in ]0,1[$.
For this, let $\delta =\bar\alpha- \alpha^\prime(0) = \bar \alpha
(1-\cos \bar \beta) \geq 0$. A simple exercise shows that the
polynomial $\alpha(s) = a s^3 + b s^2 + c s + d$ meeting the
restrictions $\alpha^\prime(0) = \alpha^\prime(1)$ and
$\alpha(1)-\alpha(0) = \bar\alpha$ is such that $a  =  - 2
\delta$, $b  =  3 \delta$, $c =  \alpha^\prime(0)$ and $d  =
\alpha(0)$. In particular $\alpha^\prime(s) = - 6 \delta s (s-1) +
\alpha^\prime(0)$. If $\cos \beta \neq 1$, then $- 6 \delta s
(s-1)
> 0$, for $s \in ]0,1[$.  As $\alpha^\prime(0) \geq 0$, then
$\alpha^\prime>0$ for $s \in ]0,1[$. If $\cos \bar\beta = 1$, then
$\delta=0$ and $\alpha^\prime(0)= \bar\alpha \cos\bar \beta >0$.
So $\alpha^\prime>0$ for $s \in [0,1]$.

\end{document}